\input amstex.tex
\input amsppt.sty
\nopagenumbers
\magnification1200

\voffset1cm
\hoffset1cm
\vsize17.5cm
\hsize11cm


\def\ad{\operatorname{ad}}

\def\bx{{\boxed{\phantom{\square}}\kern-.4pt}}

\def\CC{{\Bbb C}}
\def\ch{\operatorname{ch}}
\def\Chi{\operatorname{X}}
\def\com{\ts,\hskip-.5pt}
\def\CS{\CC\hskip-.5pt\cdot\hskip-1pt S}

\def\de{\delta}
\def\dim{\operatorname{dim}}

\def\End{\operatorname{End}\hskip1pt}
\def\EndCN{\End(\CC^N)}
\def\enddemos{{\ $\square$\enddemo}}
\def\ep{\varepsilon}
\def\Et{\tilde E}

\def\g{\frak{g}}
\def\glN{\frak{gl}_N}
\def\Gm{\Gamma_{\hskip-1pt\mu}}

\def\h{\frak{h}}

\def\id{\operatorname{id}}
\def\io{\iota}

\def\la{\lambda}
\def\La{\Lambda}
\def\ld{,\,\ldots,}

\def\mus{{\mu^{\ts\ast}}}

\def\om{\omega}
\def\ot{\otimes}

\def\ph{\varphi}

\def\Rt{\tilde R}

\def\Sg{\operatorname{S}(\g)}
\def\sgn{\operatorname{sgn}\ts}
\def\si{\sigma}
\def\soN{\frak{so}_N}
\def\spN{\frak{sp}_N}

\def\tr{\operatorname{tr}}
\def\ts{\hskip1pt}

\def\U{\operatorname{U}}
\def\Ug{\U(\g)}
\def\UglN{\U(\glN)}

\def\Z{\operatorname{Z}}
\def\Zg{\Z(\g)}
\def\ZglN{\Z(\glN)}
\def\ZZ{{\Bbb Z}}


\catcode`\@=11
\mathsurround 1.6pt

\def\hcor#1{\advance\hoffset by #1}
\def\vcor#1{\advance\voffset by #1}
\let\bls\baselineskip  \let\ignore\ignorespaces
\def\vsk#1>{\vskip#1\bls} \let\adv\advance 
\def\vv#1>{\vadjust{\vsk#1>}\ignore} \def\vvv#1>{\vadjust{\vskip#1}\ignore}
\def\vvn#1>{\vadjust{\nobreak\vsk#1>\nobreak}\ignore}
\def\vvvn#1>{\vadjust{\nobreak\vskip#1\nobreak}\ignore}
\def\setnormalbls{\edef\normalbls{\bls\the\bls}}
\def\setmaths{\edef\maths{\mathsurround\the\mathsurround}}

\def\nn#1>{\noalign{\vskip #1pt}} \def\NN#1>{\openup#1pt}
 
\let\Sum\sum \def\sum{\Sum\limits} 
\let\Prod\prod \def\prod{\Prod\limits} \let\Int\int \def\int{\Int\limits}

\let\=\m@th \def\&{.\kern.1em} \def\>{\!\;} \def\:{\!\!\;}

\ifx\plainfootnote\undefined \let\plainfootnote\footnote \fi
\expandafter\ifx\csname amsppt.sty\endcsname\relax
 
\else \fi

\newbox\s@ctb@x
\def\s@ct#1 #2\par{\removelastskip\vsk>
 \vtop{\bf\setbox\s@ctb@x\hbox{#1} \parindent\wd\s@ctb@x
 \ifdim\parindent>0pt\adv\parindent.5em\fi\item{#1}#2\strut}%
 \nointerlineskip\nobreak\vtop{\strut}\nobreak\vsk-.4>\nobreak}

\newbox\t@stb@x
\def\gadv{\global\advance} \def\gad#1{\gadv#1 1} 
\def\l@b@l#1#2{\def\n@@{\csname #2no\endcsname}%
 \if *#1\gad\n@@ \expandafter\xdef\csname @#1@#2@\endcsname{\the\Sno.\the\n@@}%
 \else\expandafter\ifx\csname @#1@#2@\endcsname\relax\gad\n@@
 \expandafter\xdef\csname @#1@#2@\endcsname{\the\Sno.\the\n@@}\fi\fi}
\def\l@bel#1#2{\l@b@l{#1}{#2}\?#1@#2?}
\def\?#1?{\csname @#1@\endcsname}
\def\[#1]{\def\n@xt@{\ifx\t@st *\def\n@xt####1{{\setbox\t@stb@x\hbox{\?#1@F?}%
 \ifnum\wd\t@stb@x=0 {\bf???}\else\?#1@F?\fi}}\else
 \def\n@xt{{\setbox\t@stb@x\hbox{\?#1@L?}\ifnum\wd\t@stb@x=0 {\bf???}\else
 \?#1@L?\fi}}\fi\n@xt}\futurelet\t@st\n@xt@}
\def\(#1){{\rm\setbox\t@stb@x\hbox{\?#1@F?}\ifnum\wd\t@stb@x=0 ({\bf???})\else
 (\?#1@F?)\fi}}
\def\dff{\expandafter\d@f} \def\d@f{\expandafter\def}
\def\edff{\expandafter\ed@f} \def\ed@f{\expandafter\edef}

\newcount\Sno \newcount\Lno \newcount\Fno
\def\Section#1{\gadno\Fno=0\Lno=0\s@ct{\the\Sno.} {#1}\par} \let\Sect\Section
\def\section#1{\gad\Sno\Fno=0\Lno=0\s@ct{} {#1}\par} \let\sect\section
\def\l@F#1{\l@bel{#1}F} \def\<#1>{\l@b@l{#1}F} \def\l@L#1{\l@bel{#1}L}
\def\Tag#1{\tag\l@F{#1}} \def\Tagg#1{\tag"\llap{\rm(\l@F{#1})}"}
\def\Th#1{Theorem \l@L{#1}} \def\Lm#1{Lemma \l@L{#1}}
\def\Prop#1{Proposition \l@L{#1}}
\def\Cr#1{Corollary \l@L{#1}} \def\Cj#1{Conjecture \l@L{#1}}
 
\def\Proof#1.{\demo{\it Proof #1}}

\def\Par{\par\medskip} \def\setparindent{\edef\Parindent{\the\parindent}}
\def\Appendix{\Sno=64\let\p@r@\z@ 
 \def\Section##1{\gad\Sno\Fno=0\Lno=0 \s@ct{} \hskip\p@r@ Appendix
\\the\Sno
  \if *##1\relax\else {.\enspace##1}\fi\par} \let\Sect\Section
 \def\section##1{\gad\Sno\Fno=0\Lno=0 \s@ct{} \hskip\p@r@ Appendix%
  \if *##1\relax\else {.\enspace##1}\fi\par} \let\sect\section
 \def\l@b@l##1##2{\def\n@@{\csname ##2no\endcsname}%
 \if *##1\gad\n@@
 \expandafter\xdef\csname @##1@##2@\endcsname{\char\the\Sno.\the\n@@}%
 \else\expandafter\ifx\csname @##1@##2@\endcsname\relax\gad\n@@
 \expandafter\xdef\csname @##1@##2@\endcsname{\char\the\Sno.\the\n@@}\fi\fi}}


\font\bigbf=cmbx10 scaled 1200
\centerline{\bigbf Capelli elements in the classical}
\smallskip
\centerline{\bigbf universal enveloping algebras} 
\bigskip\medskip
\centerline{\bf Maxim Nazarov}
\bigskip\bigskip

For any complex classical group $G=O_N,Sp_{\ts N}$ consider
the ring $\Zg$ of $G\!$-invariants
in the corresponding enveloping algebra $\Ug$.
Let $u$ be a complex parameter.
For each $n=0,1,2,\ldots$ and every partition $\nu$ of $n$  
into at most $N$ parts we define a certain rational function
$Z_\nu(u)$ which takes values in $\Zg$. 
Our definition is motivated by the works of Cherednik and Sklyanin
on the reflection equation, and also by the classical Capelli identity. 
The degrees in $\Ug$ of the values of $Z_\nu(u)$ do not exceed $n$.
We describe the images of these values in the
$n\hskip-.6pt$-th symmetric power of $\g$.
Our description involves the plethysm coefficients as studied by Littlewood,
see Theorem 3.4 and Corollary 3.6 here. 


\bigskip\section{\bf\S\ts1. Capelli elements in the algebra $\UglN$}

We work with the general linear
Lie algebra $\glN$ over the complex field $\CC$.
In this introductory section we recall the definition from
[OO1,\,S]
of the Capelli elements in the universal enveloping algebra
$\UglN$.
Here we also recall an explicit construction from [N2,O]
of these elements.

Let the indices $i,j$ run through the set $\{1\ld N\}$.
Let the vectors $e_i$ form the standard basis in $\CC^N$.
We fix in the Lie algebra $\glN$ the basis of the standard
matrix units $E_{ij}$.
We will also regard $E_{ij}$
as generators of the universal enveloping algebra $\UglN$.
Now choose the Borel subalgebra in $\glN$ spanned by the elements
$E_{ij}$ with $i\leqslant j$. Then choose the basis $E_{11}\ld E_{NN}$
in the corresponding Cartan subalgebra. 

Let $\nu$ be any partition of $n$ into at most $N$ parts.
We will write $\nu=(\nu_1\ld\nu_N)$. Let $U_\nu$ be the
irreducible $\glN$-module of highest weight $\nu$. 
The module $U_\nu$ appears in the decomposition of the 
$n$-th tensor power of the defining $\glN$-module $\CC^N$.
It is called the {\it polynomial\/} $\glN$-module corresponding
to the partition $\nu$.

There is a distinguished
basis in the centre $\ZglN$ of the universal enveloping
algebra $\UglN$, parametrized by the same partitions $\nu$. The element
$C_\nu$ of this basis is determined up to multipler from $\CC$
by the following proposition. This proposition is due to Sahi
\text{[\ts S,\,Theorem 1].}
Consider the canonical ascending filtration on the algebra $\UglN$.
With respect to this filtration
the subspace $\glN\subset\UglN$ has degree one. 

\proclaim{\quad Proposition 1.1}
There is an element $C_\nu$ in $\ZglN$ of degree at most $n$
such that for any partition $\la=(\la_1\ld\la_N)$ of not more than $n$ 
we have $C_\nu\cdot U_\la\neq\{0\}$ if and only if $\la=\nu$.
\endproclaim

We will call $C_\nu\in\ZglN$ the {\it Capelli element\/}
in the algebra $\UglN$ corresponding to the partition $\nu$.
The elements $C_\nu$ corresponding to
the partitions $\nu=(1\ld1,0\ld0)$
were studied by Capelli in [C].
In the case $\nu=(n,0\ld0)$ they were studied in 
[N1]. An explicit formula for the eigenvalue of the central element $C_\nu$
in the $\glN$-module $U_\la$ for any
$\la$ and $\nu$ was given by Okounkov and Olshanski in [\ts OO1\ts].
Let us reproduce this formula, it will fix the
multiplier from $\CC$
up to which the element $C_\nu\in\ZglN$ has been determined so far.  

Let $a=(a_1,a_2,\,\ldots)$ be an arbitrary sequence of complex numbers.
For each $k=0,1,2,\,\ldots$ introduce the $k$-th
{\it generalized factorial power}
$
(u\ts|\ts a)^k=(u-a_1)\cdots (u-a_k)
$
of the variable $u$.
Consider the function in $N$ independent variables $y_1\ld y_N$
$$
s_{\nu}(y_1\ld y_N\ts|\ts a)=
\frac
{\det\bigl[\ts(y_j|a)^{\nu_i+N-i}\ts\bigr]}
{\det\bigl[\ts(y_j|a)^{N-i}\ts\bigr]}
\Tag{1.1}
$$
where the determinants are taken with respect to $i,j=1\ld N$. This 
function is a symmetric polynomial in $y_1\ld y_N$ which is called
the {\it generalized factorial Schur polynomial,}
see [\ts M, Example I.3.20\ts].
Note that here the denominator
$$
\det\bigl[\ts(y_j|a)^{N-i}\ts\bigr]\,=\,\prod_{i<j}\,(y_i-y_j)
$$
is the Vandermonde determinant. Thus the denominator in \(1.1)
does not depend on the sequence $a$.

If $a=(0,0,\ldots)$ the polynomial
$s_\nu(y_1\ld y_N\ts|\ts a)$ is the ordinary Schur polynomial
$s_\nu(y_1\ld y_N)$.
For the general sequence $a$ by \(1.1)
$$
s_{\nu}(y_1\ld y_N\ts|\ts a)=s_\nu(y_1\ld y_N)\,+\,
\text{\rm lower degree terms}.
$$
Therefore all the polynomials $s_{\nu}(y_1\ld y_N\ts|\ts a)$
where the partitions $\nu$ have not more than $N$ parts,
form a linear basis in the ring
of symmetric polynomials in the variables $y_1\ld y_N$
with complex coefficients.

\proclaim{\quad Proposition 1.2}
The Capelli element $C_\nu\in\ZglN$ can be chosen so that its
eigenvalues in the irreducible $\glN$-modules $U_\la$ are respectively
$$
s_{\nu}\ts
(\ts\la_1+N-1\ts,\la_2+N-2\,\ld\la_N\,\,|\,\,0\ts,1\ts,2\ts,\,\ldots\,\ts)\,.
\Tag{1.2}
$$
\endproclaim

By the Harish-Chandra theorem [\ts D,\,Theorem 7.4.5\ts]\ts,
the eigenvalue of any element from $\ZglN$ in the irreducible
module $U_\la$ is a symmetric polynomial in
$\la_1+N-1\ts,\la_2+N-2\,\ld\la_N$ and all the symmetric polynomials
arise in this way. The proof of Proposition 1.2 consists of a direct
verification that when $\la_1+\ldots+\la_N\leqslant n$, 
the expression \(1.2) vanishes unless $\la=\nu$.
The details can be found in [\ts OO1, Section 3\ts]\ts.

An explicit formula for the element $C_\nu\in\UglN$ in terms
of the generators $E_{ij}$ was given by \text{[\ts N2,\, Theorem 5.3\ts]}
and \text{[O,\,Theorem 1.3].} It generalizes the
formula from [C] for $C_\nu$ with $\nu=(1\ld 1,0\ld0)$
and employs the classical results of Young [Y1,Y2] about
the irreducible representations of the symmetric group $S_n$.
Let us recall the relevant results from [Y1,Y2] here.

Let $W_\nu$ be the irreducible $S_n$-module
corresponding to the partition $\nu$.
We identify the partition $\nu$ with its Young diagram.
Fix the chain
$$
S_1\subset S_2\subset\cdots\subset S_n
\nopagebreak
\Tag{1.333}
$$
of subgroups with the standard embeddings.
There is a 
decomposition of the space $W_\nu$ into the direct sum of
one-dimensional subspaces associated with this chain. These subspaces
are parametrized by the {\it standard tableaux\/} of shape $\nu$. Each of
these tableaux is a bijective filling of the boxes of the Young diagram $\nu$
with the
numbers $1\ld n$ such that in every row and column the numbers
increase from left to right and from top to bottom
respectively. Denote by $\Cal T_\la$ the set of these tableaux.

For every tableau $T\in\Cal{T}_\la$ define a 
one-dimensional subspace $W_T$ in $W_\nu$ as follows. For any
$p\in\{1\ld n\}$ take the tableau obtained from $T$ by
removing each of the numbers $p+1\ld n$. Let the Young diagram $\om$
be its shape. The subspace $W_T$ is contained in an irreducible
$S_p$-submodule of $W_\nu$ corresponding to $\om$. Any basis of
$W_\nu$ formed by vectors $w_T\in W_T$ is called a 
{\it Young basis}. Fix an $S_n$-invariant inner product
$\langle\,\,,\,\rangle_{\nu}$ 
in $W_\nu$. The subspaces
$W_T$ are then pairwise orthogonal. We shall be assuming that
$\langle\ts w_T,w_T\ts\rangle_\nu=1$ for each tableau $T\in\Cal{T}_\la$.

For any tableau $T\in\Cal{T}_\nu$
consider the normalized diagonal matrix element of the
$S_n$-module $W_\nu$ corresponding to the vector $w_T$
$$
\Phi_T\,=\,\frac{\dim W_\nu}{n!}
\sum_{\si\in S_n}\,
\langle\,w_T\com\ts\si\!\cdot\!\hskip.5pt w_T\,\rangle_\nu
\,\si
\,\in\,\CS_n\,.
\Tag{1.0}
$$
This formula is the most simple when $T\in\Cal{T}_\nu$ is
the {\it column tableau\/}. This tableau is obtained by
filling the boxes of the diagram $\nu$ with $1\ld n$ by columns from left
to right, downwards in each column. We shall denote
this tableau by $T_c$.
Let $S_\nu$ and $S_\nu^{\ts\prime}$ be the subgroups
in $S_n$ preserving the collections of numbers appearing respectively
in every row and column of the tableau $T_c$.
Take the elements of the group ring $\CS_n$
$$
\Theta_\nu=\sum_{\si\in S_\nu}\si
\ \quad\text{and}\ \quad
\Theta_\nu^{\ts\prime}=\sum_{\si\in S_\nu^{\ts\prime}}\si\cdot\sgn\si\,.
$$
As usual, we denote by $\nu_1^{\ts\prime},\nu_2^{\ts\prime}\ts,\,\ldots$
the column lenghts \text{of the diagram $\nu$.} Then by [Y1]
$$
\Phi_{T_c}=\,
\frac{\dim W_\nu}{n!}\cdot
\frac{\Theta_\nu^{\ts\prime}\Theta_\nu\Theta_\nu^{\ts\prime}}
{\nu_1^{\ts\prime}\ts!\,\nu_2^{\ts\prime}\ts!\,\ldots}\ .
$$

There is an alternative description of the one-dimensional subspace
$W_T$ in $W_\nu$ due to Jucys [\ts J\ts].
Consider the sum of transpositions 
$$
z_p=(1\com p)+(2\com p)+\ldots+(\ts p-1\com p)\,\in\,\CS_n\,.
$$
The elements $z_1\ld z_n\in\CS_n$ are called  the {\it Jucys-Murphy elements}
corresponding to the standard chain \(1.333). They pairwise commute.
Fix a tableau $T\in\Cal{T}_\nu$. For every $r=1\ld n$
put $c_p=k-l$ if the number $p$ appears in the
$k$-th column and $l$-th row of the tableau $T$. The number
$c_p$ is called the {\it content} of the box of the diagram $\nu$
occupied by $p$. Here on the left we show
the column tableau of shape $\nu=(4,3,1)$:
\smallskip\smallskip\smallskip\smallskip
\vbox{
$$
{\bx}
{\bx}
{\bx}
{\bx}
\phantom{\bx}
\phantom{\bx}
\phantom{\bx}
\phantom{\bx}
{\bx}
{\bx}
{\bx}
{\bx}
$$
\vglue-17.8pt
$$
{\bx}
{\bx}
{\bx}
\phantom{\bx}
\phantom{\bx}
\phantom{\bx}
\phantom{\bx}
\phantom{\bx}
{\bx}
{\bx}
{\bx}
\phantom{\bx}
$$
\vglue-17.9pt
$$
{\bx}
\phantom{\bx}
\phantom{\bx}
\phantom{\bx}
\phantom{\bx}
\phantom{\bx}
\phantom{\bx}
\phantom{\bx}
{\bx}
\phantom{\bx}
\phantom{\bx}
\phantom{\bx}
$$
\vglue-57.5pt
$$
1
\kern9pt
4
\kern9pt
6
\kern9pt
8
\kern67pt
0
\kern9pt
1
\kern9pt
2
\kern9pt
3
$$
\vglue-18pt
$$
2
\kern9pt
5
\kern9pt
7
\kern3pt
\phantom8
\kern70pt
\text{-1}
\kern9pt
0
\kern9pt
1
\kern9pt
\phantom2
$$
\vglue-18pt
$$
3
\kern9pt
\phantom3
\kern9pt
\phantom6
\kern9pt
\phantom8
\kern64pt
\text{-2}
\kern9pt
\phantom0
\kern9pt
\phantom0
\kern9pt
\phantom0
$$
}
\smallskip\smallskip\smallskip
\noindent
On the right we have indicated the contents of the boxes of the
Young \text{diagram $\nu=(4,3,1)$.}
So here we get
$
(c_1,\ldots\hskip0pt,c_8)\hskip-1pt=\hskip-1pt
\text{(0\ts,\ts-1\ts,\ts-2\ts,\ts1\ts,\ts0\ts,\ts2\ts,\ts1\ts,\ts3)}\!$.
Observe that the standard
tableau $T\in\Cal{T}_\nu$ can be always recovered from the sequence
of contents $c_1\ld c_n$. The next lemma is \text{contained in [\ts J\ts].}

\proclaim{\quad Lemma 1.3}
We have $z_p\cdot w_T=c_p\ts w_T$ in $W_\nu$ for any $p=1\ld n$.
\endproclaim

Let us now reproduce the explicit formula from [N2,O] for the
element $C_\nu\in\UglN$. Consider the permutational action
of the symmetric group $S_n$ in the tensor product $(\CC^N)^{\ot n}$.
Denote by $Y_T$ the linear operator in $(\CC^N)^{\ot n}$ 
corresponding of the element \(1.0)\ts.
The image of this operator is equivalent to $U_\nu$ as a $\glN$-module,
see \text{[\ts W, Section IV\hskip-1pt.4\ts]\ts.} 
Moreover, by definition we have the equality $Y_T^{\ts2}=Y_T$. 

Further, denote by $\io_p$ the embedding of the algebra $\EndCN$ into the
tensor product $\EndCN^{\ot n}$ \text{as the $p$-th tensor factor:}
$$
\io_p(X)=1^{\ot\ts(p-1)}\ot X\ot1^{\ot\ts(n-p)}\ts;
\qquad
p=1\ld n\ts.
\Tag{1.225}
$$
We will use this notation throughout the present article. Now put 
$$
\align
E(u)\,&=\,-\,u\,+\,\sum_{ij}\,E_{ij}\ot E_{ji}\,\in\,
\EndCN\ot\UglN\ts[\ts u\ts]\,,
\hskip-100pt
\Tag{1.24}
\\
E_p(u)\,&=\,(\io_p\ot\id)\,\bigl(\ts E(u)\bigr)\in
\EndCN^{\ot n}\ot\Ug\,[\ts u\ts]\ts.
\Tag{1.25}
\endalign
$$
Let $\tr:\EndCN\to\CC$ be the usual matrix trace,
so that $\tr\ts(E_{ij}$) equals the Kronecker delta $\de_{ij}$. 
Now consider the product
$$
E_1(u_1)\cdots E_n(u_n)\,\in\,\EndCN^{\ot n}\ot\UglN[u_1\ld u_n]\,.
$$

\proclaim{\quad Theorem 1.4}
For any standard tableau $T\in\Cal{T}_\nu$
we have
$$
C_\nu=(\tr^{\ts\ot n}\ot\id)
\bigl(\ts Y_T\ot1\cdot E_1(c_1)\cdots E_n(c_n)\ts\bigr)\,.
\Tag{1.3}
$$
\endproclaim

The proofs of this theorem given in [N2,O] were rather involved.
More elegant proof was subsequently found by Molev
\text{[\ts M2, Theorem 8.2\ts].} 
All these results were based on the notion of a fusion procedure
\text{introduced} by Cherednik in [C2]. We keep using this notion
in the present article.

Consider again the group ring $\CS_n$. For every two distinct indices
$p\ts,q=1\ld n$ introduce the rational
function of two complex variables $u\ts,v$ valued in $\CS_n$
$$
\ph_{pq}(u\com v)=1-\frac{(\ts p\,,q\ts)}{u-v}\,.\hskip-10pt
$$
As direct calculations show, these rational function satisfy the equations
$$
\ph_{pq}(u\com v)\,\ph_{pr}(u\com w)\,\ph_{qr}(v\com w)=
\ph_{qr}(v\com w)\,\ph_{pr}(u\com w)\,\ph_{pq}(u\com v)
\Tag{1.4}
$$
for all pairwise distinct indices $p\com q\com r$.
Consider the rational function of $u\com v\com w$ appearing at either side
of \(1.4)\ts.
The factor
$\ph_{pr}(u\com w)$ in \(1.4) has a pole at $u=w$. However, we have the
following lemma.

\proclaim{\quad Lemma 1.5}
The restriction of \(1.4) to the set of all $(u\com v\com w)$
such that $v=w\pm1$, is regular at $u=w$.
\endproclaim

\demo{\quad Proof}
Under the condition $v=w\pm1$ the rational function \(1.4) can be
written as 
$$
\biggl(1-\frac{(p\com q)+(p\com r)}{u-w\mp1}\biggr)\cdot
\bigl(\ts 1\mp(q\com r)\bigr)
$$
which is a rational function of $u\com w$ manifestly regular at $u=w$
\enddemos

Using Lemma 1.5 one can prove the next proposition,
for details see [\ts N2\ts,\ts Proposition 2.12\ts].
Let the superscript ${}^\vee$ denote the
group embedding $S_n\to S_{n+1}$ determined by
the assignment 
$(p\com q)\mapsto(\ts p+1\com q+1)$.

\proclaim{\quad Proposition 1.6}
We have the identity in the algebra $\CS_{n+1}(u)$
$$
\Bigl(\ts1-\sum_{p=1}^n\,\frac{(1,p+1)}{u}\ts\Bigr)
\cdot\Phi_T^\vee\ts=
\ph_{12}(u\com c_1)\ts\cdots\ts\ph_{1,n+1}(u\com c_n)
\cdot\Phi_T^\vee\,.
$$
\endproclaim

The proof of the next proposition is similar and will be also omitted.

\proclaim{\quad Proposition 1.7}
We have the identity in the algebra $\CS_{n+1}(u)$
$$
\Bigl(\ts1+\sum_{p=1}^n\,\frac{(p\ts,n+1)}{u}\ts\Bigr)
\cdot\Phi_T\ts=
\ph_{1,n+1}(-c_1\com u)\ts\cdots\ts\ph_{n,n+1}(-c_n\com u)\cdot\Phi_T\,.
\hskip-10pt
$$
\endproclaim

We will also use an alternative definition of the
element $\Phi_T\in\CS_n$ due to Cherednik [C2].
Suppose the numbers $1\ld n$ appear respectively in the rows $l_1\ld l_n$
of the standard tableau $T$. Order the set of all pairs $p,q$
with $1\leqslant p<q\leqslant n$ lexicographically.

\proclaim{\quad Theorem 1.8}
The rational function of $u\!$ defined as the ordered product in $\CS_n(u)$
of the elements
$\ph_{pq}(\ts c_p+l_p\ts u\com c_q+l_q\ts u\ts)$
over the pairs $p\com q$
is regular at $u=0$, and takes at $u=0$ the value $\Phi_T$.
\endproclaim

One can prove this theorem
by again using Lemma 1.5. This proof is contained in [N2,\,Section 2].
Another proof can be found in [JKMO].

We will close this section with a generalization
of Theorem~1.4. Let us consider 
for a standard tableau $T\in\Cal T_\nu$ the element of $\UglN[u]$
$$
(\tr^{\ts\ot n}\ot\id)
\bigl(\ts Y_T\ot1\cdot E_1(u+c_1)\cdots E_n(u+c_n)\bigr)\,.
\Tag{1.11}
$$

\proclaim{\quad Corollary 1.9}
The element \(1.11) belongs to $\ZglN[u]$ and does not depend on the choice
of a tableau $T\in\Cal T_\nu$. The eigenvalue of \(1.11) in the irreducible
$\glN$-module $U_\lambda$ is
$$
s_{\nu}\ts
(\ts\la_1-u+N-1\ts,\la_2-u+N-2\,\ld\la_N-u\,\,|
\,\,0\ts,1\ts,2\ts,\,\ldots\,\ts)\,.
$$
\endproclaim

\vbox{
\demo{\quad Proof}
For any complex value of the parameter $u$ consider the automorphism
of the unital algebra $\UglN$ determined by the assignment
$E_{ij}\mapsto E_{ij}-u\cdot\de_{ij}$.
The element \(1.11) can be obtained by applying this
automorphism to the left hand side of \(1.3), see the definition \(1.24).
So the first statement of Corollary 1.9 follows from Theorem 1.4.
By pulling back the $\glN$-module $U_\la$ through that automorphism we obtain
the irreducible $\glN$-module of highest weight 
$(\la_1-u\ld\la_N-u)$. The second statement of Corollary 1.9
now follows from Proposition 1.2
\enddemos

The principal aim of this article is to introduce the analogues of the
elements \(1.11) for the remaining classical Lie algebras $\soN$ and $\spN$.
}


\bigskip\section{\bf\S\ts2. Traceless tensors in the space $(\CC^N)^{\ot m}$}

We will regard the orthogonal and symplectic
Lie algebras $\soN$ and $\spN$ as subalgebras in $\glN$.
{}From now on we will let the indices $i,j$ run through the set
$\{\ts-\ts M\ld-1,1\ld M\ts\}$ if $N=2M$ and through the set
$\{\ts-\ts M\ld-1,0,1\ld M\ts\}$ if $N=2M+1$.
Let $e_i$ be the elements of the standard basis in $\CC^N$.
We will realize the complex orthogonal group $O_N$ as the subgroup in $GL_N$
preserving
the symmetric bilinear form
$\langle\,e_i\ts,e_j\,\rangle=\de_{i,-j}$
on $\CC^N$.
The complex symplectic group $Sp_N$ 
will be realized as the subgroup in $GL_N$ preserving the
alternating form
$\langle\,e_i\ts,e_j\,\rangle=\de_{i,-j}\cdot\sgn i\,$.

Let $G$ be any of the subgroups $O_N,Sp_N$ in $GL_N$. Denote by
$\g$ the corresponding Lie subalgebra in $\glN$.
Put $\ep_{ij}=\sgn i\cdot\sgn j$ if $G=Sp_N$ and
$\ep_{ij}=1$ if $G=O_N$.
The Lie subalgebra $\g\subset\glN$ is then spanned by the elements
$$
F_{ij}=E_{ij}-\ep_{ij}\cdot E_{-j,-i}\,.
$$
We will also regard $F_{ij}$ as generators of the universal 
enveloping algebra $\Ug$.
We choose the Borel subalgebra in $\g$ spanned by the elements
$E_{ij}$ with $i\leqslant j$.
Let us fix  the basis $F_{-M,-M}\ld F_{-1,-1}$
in the corresponding Cartan subalgebra $\h\subset\g$. 
Any weight $\mu=(\mu_1\ld\mu_M)$ of $\h$ will be taken with respect to this
basis. The half-sum of the positive roots of $\h$ is
$$
\rho=(\ep+M-1,\ep+M-2\ld\ep)
$$
where $\ep=0\,,\,\frac12\,,\,1$
for $\g=\frak{so}_{2M}\,,\,\frak{so}_{2M+1}\,,\,\frak{sp}_{2M}$ respectively.

Now we assume that $\mu$ is a partition of $m$ with at most $M$ parts.
Then $\mu$ can be regarded as a dominant weight of $\h$. Let $V_\mu$
be the irreducible $\g$-module of the highest weight $\mu$.
Note that if $\g=\frak{so}_{2M}$ then
$\mus=(\mu_1\ld\mu_{M-1},-\ts\mu_M)$
is again a dominant weight of $\h$. We will also consider
the corresponding irreducible $\frak{so}_{2M}$-module $V_\mus$.
It is obtained from the module $V_\mu$
via the conjugation in $\frak{so}_{2M}\subset\EndCN$ by
$$
E_{1,-1}+E_{-1,1}+E_{22}+E_{-2,-2}+\ldots+E_{MM}+E_{-M,-M}\,\in\,O_N\,.
\hskip-28pt
$$

All these irreducuble $\g$-modules appear in the decomposition of the 
$m$-th tensor power of the identity $\g$-module $\CC^N$.
Take any two distinct indices $p$ and $q$ from the set $\{1\ts\ld m\}$.
By applying the $G\!$-invariant bilinear form
$\langle\,\,,\,\rangle$ on $\CC^N$ to
an element $t\in(\CC^N)^{\ot m}$ in the $p$-th and \text{$q$-th} tensor factors
we obtain a certain element $t^{\ts\prime}\in(\CC^N)^{\ot(m-2)}$.
Then the element $t$ is called {\it traceless} if $t^{\ts\prime}=0$ for all
possible indices $p\neq q$.

Now fix any embedding of the irreducible $\glN$-module
$U_\mu\!$ to $(\CC^N)^{\ot m}\!$. The subspace $U_\mu\cap V$
in $U_\mu$ is preserved by the action of the subalgebra $\g\subset\glN$.
For $\g=\frak{so}_{2M+1}\ts,\ts\frak{sp}_{2M}$ this subspace is isomorphic
to $V_\mu$ as $\g$-module. For $\g=\frak{so}_{2M}$ it is isomorphic to $V_\mu$
only if $\mu_M=0$. Otherwise $U_\mu\cap V$ splits into the direct sum of the
$\frak{so}_{2M}$-modules $V_\mu$ and $V_\mus$.
All these statements are contained in
\text{[\ts W,\,Section V\hskip-0.4pt.\ts9\ts]\ts.}

We denote by $\Zg$ the ring of invariants in the universal
enveloping algebra $\Ug$ with respect to the adjoint action of the group $G$.
The ring $\Zg$ coincides with the centre of $\Ug$ for
$\g=\frak{so}_{2M+1}\ts,\ts\frak{sp}_{2M}$ but is strictly contained
in the centre of $\Ug$ when $\g=\frak{so}_{2M}$.
Then any element of $\Z(\frak{so}_{2M})$ acts in the irreducible
$\frak{so}_{2M}$-modules $V_\mu$ and $V_\mus$ by the same scalars. 

There is a distinguished basis in the vector space $\Zg$ analogous
to the basis of the Capelli elements $C_\nu$ in $\ZglN$.
This basis is labelled by the partitions $\mu$ and was introduced
by Okounkov and Olshanski by generalizing Proposition 1.1.  
The element $B_\mu$ of this basis is determined up to multipler from $\CC$
by the next proposition [OO2, \!Theorem 2.3\ts]\ts. 
Consider the canonical ascending filtration on the algebra $\Ug$.
With respect to this filtration
the subspace $\g\subset\Ug$ has degree one. 

\proclaim{\quad Proposition 2.1}
There exists an element $B_\mu$ in $\Zg$ of degree at most $2\ts m$
such that for any partition $\la=(\la_1\ld\la_M)$ of not more than $m$ 
we have $B_\mu\cdot V_\la\neq\{0\}$ if and only if $\la=\mu$.
\endproclaim

Explicit formula for the eigenvalue of the element $B_\mu\in\Zg$
in the irreducible $\g$-module $V_\la$ for any
$\la$ and $\mu$ has been also given in [OO2].
We will reproduce this formula, it fixes the multiplier from $\CC$
up to which the element $B_\mu\in\Zg$ is determined by Proposition 2.1.
This formula again employs the definition \(1.1).  

\proclaim{\quad Proposition 2.2}
The element $B_\mu\in\Zg$ can be chosen so that its
eigenvalues in the irreducible $\g$-modules $V_\la$ are respectively
$$
s_{\mu}\ts
\bigl(\ts(\la_1+\rho_1)^2\ld(\la_M+\rho_M)^2
\,\,|\,\,\ep^2,(\ep+1)^2,\,\ldots\,\ts\bigr)\,.
$$
\endproclaim

The proof of this proposition does not differ significantly
from that of Proposition 1.2. For details see [OO2, Theorem 2.5].
A certain explicit expression for the element $B_\mu\in\Ug$
in terms of the generators $F_{ij}$
has been recently given by Olshanski in [\ts O2\ts ].
This is an analogue of the expression [OO1\,,\,Theorem 14.1]
for the element
$C_\nu\in\UglN$ which is more complicated than \(1.3).
An analogue of the formula \(1.3) for $B_\mu$ with the 
general partition $\mu$ is unknown.     
For $\mu=(1\ld 1,0\ts\ld0)$ and $\mu=(m\ts,0\ts\ld0)$
this analogue was given in [MN]. In the present article we will
consider a natural generalization of the construction [MN].
But in general it yields elements of the ring $\Zg$ different from $B_\mu$. 

Similarly to \(1.225), for any element $X\in\EndCN^{\ot2}$
and any two distinct indices $p\com q\in\{1\ld n\}$ with fixed $n$
we will denote
$$
X_{pq}=(\io_p\ot\io_q)\,(X)\ts\in\ts\EndCN^{\ot n}\ts.
$$
Along with Lemma 1.5, we will use one more simple observation. Denote
$$
F(u)\,=\,-\,u\,-\,\eta\,+\sum_{ij}\,E_{ij}\ot F_{ji}\,\in\,
\EndCN\ot\UglN\ts[\ts u\ts]
\Tag{2.0}
$$
where we set $\eta=\frac12\,,-\frac12$ for $\g=\soN,\spN$ respectively. Let
$$
\Et(u)\,=\,-\,u\,+\sum_{ij}\,\ep_{ij}\cdot E_{ij}\ot E_{-i,-j}
\,\in\,\EndCN\ot\UglN\ts[\ts u\ts]
$$
be the element
obtained from $E(u)$ by applying the transposition with respect
to the bilinear $\langle\,\,,\,\rangle$ in the tensor factor $\EndCN$.
Now consider
$$
\frac{\Et(\eta-u)\,E(\eta+u)}{u-\eta}\,\in\,\EndCN\ot\UglN\ts(\ts u\ts)\,.
\Tag{2.3}
$$
We have the standard representation $\UglN\to\EndCN^{\ot m}$ which
makes the element \(2.3) acting in the space $(\CC^N)^{\ot(m+1)}$.
The element $F(u)\in\EndCN\ot\Ug\ts[\ts u\ts]$ also acts in the
space $(\CC^N)^{\ot(m+1)}$ and the latter action preserves the
subspace $\CC^N\ot V$. Here is a simple lemma.

\proclaim{\quad Lemma 2.3}
\!\!The action of the element \(2.3) in the space $(\CC^N)^{\ot(m+1)}$
preserves the subspace $\CC^N\ot V$. The action of 
the element $F(u)$ in this subspace coincides with the action of\/ \(2.3).
\endproclaim

\demo{\quad Proof}
Consider the elements of the algebra $\EndCN^{\ot2}$
$$
P=\sum_{ij}\,E_{ij}\ot E_{ji}
\qquad\text{and}\qquad
Q=\sum_{ij}\,\ep_{ij}\cdot E_{ij}\ot E_{-i,-j}\,.
$$
The element $P$ corresponds to the exchange operator
$e_i\ot e_j\mapsto e_j\ot e_i$ in $(\CC^N)^{\ot2}$.
The element $Q$ is obtained from $P$ by applying to either tensor factor
of $\EndCN^{\ot2}$
transposition with respect to $\langle\,\,,\,\rangle$.
\text{Observe that}
$$
P\,Q=Q\ts P=
\cases
\phantom{-\,}Q&\text{if\ \ $\g=\soN$\ts,}
\\
-\,Q&\text{if\ \ $\g=\spN$\ts.}
\endcases
\Tag{2.35}
$$
Further, by the definition of a traceless tensor $t\in\EndCN^{\ot m}$
we have the equality $Q_{pq}\ts t=0$ for any two distinct
indices $p\com q\in\{1\ts\ld m\}$.

By definition the image of the element $F(u)\in\EndCN\ot\Ug\ts[\ts u\ts]$
in $\EndCN^{\ot(m+1)}\ts[\ts u\ts]$
under the representation \text{$\Ug\to\EndCN^{\ot m}$ is}
$$
P_{12}+\ldots+P_{1,m+1}-Q_{12}-\ldots-Q_{1,m+1}-u-\eta\,.
\Tag{2.4}
$$
On the other hand, the image of the element \(2.3) in
$\EndCN^{\ot(m+1)}\ts(\ts u\ts)$ under the
representation $\UglN\to\EndCN^{\ot\ts m}$ is the product
$$
\biggl(1+\frac{Q_{12}+\ldots+Q_{1,m+1}}{u-\eta}\biggr)
\cdot
\bigl(P_{12}+\ldots+P_{1,m+1}-u-\eta\ts\bigr)\ts.
\nopagebreak
$$
By \(2.35) and the definition of $\eta$ 
this product equals \(2.4) plus the sum
$$
\sum_{p\neq q}\ \,\frac{Q_{1,q+1}\,P_{1,p+1}}{u-\eta}
\,=\,
\sum_{p\neq q}\ \,\frac{P_{1,p+1}\,Q_{p+1,q+1}}{u-\eta}\,\,.
\nopagebreak
$$
But the action of the latter sum in $\CC^N\ot V$ is identically zero
\enddemos

Using this lemma we can easily prove the following proposition.
It is a particular case of a more general result from [O1]. Denote
$$
R(u\com v)\ts=\ts1-\frac{P}{u-v}
\qquad\text{and}\qquad\!
\Rt(u\com v)\ts=\ts1+\frac{Q}{u+v}
$$
in $\EndCN^{\ot2}(u,v)$. The first of these two functions
is the {\it rational Yang $R$-matrix\/}.
For any two distinct indices $p\com q\in\{1\ld n\}$
the element $R_{pq}(u\com v)\in\EndCN^{\ot n}(u\com v)$ corresponds
to $\ph_{pq}(u\com v)\in\CS_n(u\com v)$
under the permutational action of the symmetric group $S_n$ in
$(\CC)^{\ot n}$.

Similarly to \(1.25), for any fixed $n$ and every index $p=1\ld n$
let $F_p(u)\in\EndCN^{\ot n}\ot\Ug[\ts u\ts]$ and 
$\Et_p(u)\in\EndCN^{\ot n}\ot\UglN[\ts u\ts]$
be the images
of $F(u)$ and $\Et(u)$ with respect to the embedding $\io_p\ot\id$.
In the equations 
(2.5) 
to
(2.9) 
below we will write $R(u,v)$ and $\Rt(u,v)$
instead of $R(u,v)\ot1$ and $\Rt(u,v)\ot1$ in $\EndCN^{\ot2}\ot\Ug$ for short.

\proclaim{\quad Proposition 2.4}
We have the relation in\/ $\EndCN^{\ot2}\ot\Ug\,(u,v)$
$$
R(u,v)\,F_1(u)\,\Rt(u,v)\,F_2(v)=F_2(v)\,\Rt(u,v)\,F_1(u)\,R(u,v)\,.
\Tag{2.5}
$$
\endproclaim

\demo{\quad Proof}
This proposition can be verified by direct calculation. However, 
we give a conceptual proof which comes back to the origin [C2\ts,\ts S2] of the
{\it reflection equation\/} \(2.5). 
In the algebra $\EndCN^{\ot2}\ot\UglN\,(u,v)$ 
$$
\align
R(u,v)\,E_1(u)\,E_2(v)&\,=\,E_2(v)\,E_1(u)\,R(u,v)\,,
\Tag{2.81}
\\
R(u,v)\,\Et_1(-u)\,\Et_2(-v)&\,=\,\Et_2(-\,v)\,\Et_1(-u)\,R(u,v)\,,
\Tag{2.82}
\\
\Et_1(-u)\,\Rt(u,v)\,E_2(v)&\,=\,E_2(v)\,\Rt(u,v)\,\Et_1(-u)\,,
\Tag{2.83}
\\
E_1(u)\,\Rt(u,v)\,\Et_2(-v)&\,=\,\Et_2(-v)\,\Rt(u,v)\,E_1(u)\,.
\Tag{2.84}
\endalign
$$
The relation \(2.81) is well known and can be easily
verified. The relation \(2.82) is obtained from \(2.81) by applying
in the tensor factor $\UglN$ the automorphism 
$E_{ij}\mapsto-\ts\ep_{ij}\cdot E_{-j,-i}$.
Applying to \(2.81) transposition with respect to
$\langle\,\,,\,\rangle$ in the first tensor factor of $\EndCN^{\ot2}$
we obtain \(2.83).
By applying to \(2.82) transposition with respect to
$\langle\,\,,\,\rangle$ in the second tensor factor of $\EndCN^{\ot2}$
we obtain \(2.84).

Using \(2.81) to \(2.84) we get the equality
in $\EndCN^{\ot2}\ot\UglN\,(u,v)$
$$
\align
&R(u,v)\,
\Et_1(\eta-u)(u)\,E_1(\eta+u)\,
\Rt(u,v)\,
\Et_2(\eta-v)\,E_2(\eta+v)=
\\
&\Et_2(\eta-v)\,E_2(\eta+v)\,
\Rt(u,v)\,
\Et_1(\eta-u)(u)\,E_1(\eta+u)\,
R(u,v)\,\ts.
\endalign
$$
The intersection of the kernels of all the representations
$\Ug\mapsto\End V$ for $n=1,2,\ldots$ is zero
[\ts D,\, Theorem 2.5.7\ts]\ts,
therefore \(2.5) follows from the above equality
in $\EndCN^{\ot2}\ot\UglN\,(u,v)$ by Lemma 2.3
\enddemos

We will now introduce the main object of our study in this article.
Let $\nu$ be any partition of $n$ with at most
$N$ parts. Let $T$ be any standard tableau of shape $\nu$.
It determines the sequence of contents $c_1\ld c_n$. 
Consider the element of the algebra $\EndCN^{\ot n}\ot\Ug(u)$
$$
F_T(u)=
(Y_T\ot1)\ts\cdot\ts\prod_{p=1}^n\, 
\Bigl(\ts1+\frac{Q_{1p}\ot1+\ldots+Q_{p-1,p}\ot1}{2u+c_p}\ts\Bigr)
\ts F_p(u+c_p)
$$
where the (noncommuting) factors corresponding to $s=1\ld n$
are arranged from the left to right.
For example, for each of the partitions
$\nu=(2)$ and $\nu=(1,\!1)$ there is only one standard tableau
of shape $\nu$. For these partitions we get the elements 
of the algebra $\EndCN^{\ot2}\ot\Ug(u)$
$$
F_T(u)=
\bigl(\ts1\pm P\ot1\ts\bigr)\cdot F_1(u)\cdot
\Bigl(\ts1+\frac{Q\ot1}{2u\pm1}\ts\Bigr)\cdot F_2(u\pm1)
$$
respectively. Our main object of study is the rational function of $u$
$$
Z_\nu(u)=(\tr^{\ts\ot n}\ot\id)\ts\bigl(F_T(u)\bigr)
\Tag{2.6}
$$
which by definition takes values in $\Ug$, cf.\ \(1.11).
As we will show later, this function
does not depend on the choice of the tableau $T\in\Cal{T}_\nu$.

\proclaim{\quad Proposition 2.5}
The function $Z_\nu(u)$ takes values in the ring $\Zg$.
\endproclaim

\demo{\quad Proof}
We regard the group $G$ as a subgroup in $GL_N\subset\EndCN$. Consider
the adjoint action $\ad$ of the group $G$ in the enveloping algebra $\Ug$.
Observe that by the definition \(2.0) for any element $g\in G$
$$
(\id\ot\ad g)\,\bigl(F(u)\bigr)=g\ot1\cdot F(u)\cdot g^{-1}\ot1\,.
$$
Elements $Y_T\ts,Q_{1p}\ld Q_{p-1,p}\in\EndCN^{\ot n}$
commute with $g^{\ts\ot n}$. So
$$
(\id\ot\ad g)\,\bigl(F_T(u)\bigr)=
\bigl(g^{\ts\ot n}\ot1\bigr)\cdot
F_T(u)\cdot\bigl((g^{-1})^{\ot n}\ot1\bigr)\,.
\hskip-20pt
$$
Hence
$$
(\tr^{\ts\ot n}\ot\ad g)\,\bigl(F_T(u)\bigr)=
(\tr^{\ts\ot n}\ot\id)\,\bigl(F_T(u)\bigr)
\quad\square
$$
\enddemo

We need one more formula for the element
$F_T(u)$. It has motivated our definition of $Z_\nu(u)$.
We will keep
to the convention used in the definition of $F_T(u)$:
in any product over a certain index the noncommuting factors are arranged
from the left to the right, as this index increases.

\proclaim{\quad Proposition 2.6}
\!\!The element $F_T(u)\in\EndCN^{\ot n}\ot\Ug(u)$ equals
$$
(\ts Y_T\ot1\ts )\ts\cdot\prod_{p=1}^n\ 
\biggl(\ 
\prod_{q=1}^{p-1}\,
\Rt_{qp}(u+c_q,u+c_p)\ot1
\biggr)
\ts F_p(u+c_p)\ts.
$$
\endproclaim

\demo{\quad Proof}
We use the induction on $n$. In the case $n=1$ the required
equality is tautological. Assume we have the required equality for
some partition $\nu$ of $n\geqslant1$.
Take any standard tableau $U$ with $n+1$ boxes and not more
than $N$ rows, such that by removing the box with number $n+1$ we get $T$.
Let $c$ be the content of the removed box.
Consider the projector
$Y_U\in\EndCN^{\ot(n+1)}$. It is divisible on the right
by $Y_T\ot\id$. So by definition the element
$F_T(u)\in\EndCN^{\ot(n+1)}\ot\Ug(u)$ equals
$$
(Y_U\ot1)\cdot
F_T(u)\cdot
\Bigl(\ts1+\frac{Q_{1,n+1}\ot1+\ldots+Q_{n,n+1}\ot1}{2u+c}\ts\Bigr)
\cdot F_{n+1}(u+c)
$$
But the element $F_T(u)\in\EndCN^{\ot n}\ot\Ug(u)$
is divisible on the right by \text{$Y_T\ot1$.}
This follows from Theorem 1.8 and Proposition 2.4\ts,
see also [MNO, Section 4.2].
The alternative expression for $F_U(u)$ is now provided
by inductive assumption and by the identity in
$\EndCN^{\ot(n+1)}(v)$
$$
\align
(Y_T\ot\id)
&\cdot
\bigl(\,1+(\ts Q_{1,n+1}+\ldots+Q_{n,n+1})/v\,\bigr)=
\Tag{2.9}
\\
(Y_T\ot\id)
&\cdot
\Rt_{1,n+1}(c_1,v)\ts\cdots\ts\Rt_{n,n+1}(c_n,v)
\endalign
$$
with $v=2u+c$. Let us verify this identity.  
By applying to both sides of the equation \(2.9)
the transposition with respect
to $\langle\,\,,\,\rangle$ in each of the first $n$ tensor factors
in $\EndCN^{\ot(n+1)}$ we get
$$
\align
\bigl(\,1+(\ts P_{1,n+1}+\ldots+P_{n,n+1})/v\,\bigr)
&\cdot
(Y_T\ot\id)=
\Tag{2.10}
\\
R_{1,n+1}(-c_1,v)\ts\cdots\ts R_{n,n+1}(-c_n,v)
&\cdot
(Y_T\ot\id)\,.
\endalign
$$
We used the fact that the element $Y_T\ot\id\in\EndCN^{\ot(n+1)}$
is invariant under this transposition. But \(2.10) is provided
by Proposition 1.7
\enddemos

\proclaim{\quad Theorem 2.7}
Here\/ $Z_\nu(u)$ does not depend on the choice of 
$T\in\Cal{T}_\nu$.
\endproclaim

\demo{\quad Proof}
Any standard tableau of the shape $\nu$ can be obtained from
the \text{column}
tableau $T_c$ by a chain of transformations $T\mapsto T^{\ts\prime}$
where the entries of the tableaux $T,T^{\ts\prime}\in\Cal{T}_\nu$
differ by a single transposition \text{$(r\com r+1)$}
such that $l_r>l_{r+1}$ for the tableau $T$. 
Let $c_1^{\,\prime}\ld c_n^{\,\prime}$ be the sequence of contents
of the tableau $T^{\ts\prime}$. It is obtained from the sequence
$c_1\ld c_n$ by exchanging the terms $c_r$ and $c_{r+1}$.
Note that here $|\ts c_r-c_{r+1}|>1$, put $d=(\ts c_r-c_{r+1})^{-1}$.
Due to \text{[\ts Y2\ts,\, Theorem IV\ts]} we have the relation
$$
\Phi_{T^\prime}=
\bigl((\ts r,r+1)+d\,\bigr)\,
\frac{\Phi_T}{1-d^{\ts2}}\,
\bigl((\ts r,r+1)+d\,\bigr)
\Tag{2.11}
$$
in the group ring $\CS_n$, see the definition \(1.0)\,.
Let $X$ be the product $PR\ts(c_{r+1}\com c_r)=P+d$ in $\EndCN^{\ot2}$, then
the relation \(2.11) implies 
$$
Y_{T^\prime}=X_{r,r+1}\,\frac{Y_T}{1-d^{\ts2}}\ts\,X_{r,r+1}
\Tag{2.12}
$$
in $\EndCN^{\ot n}$. On the other hand, 
by using \text{Proposition 2.4 we obtain}
$$
\align
X_{r,r+1}\,\cdot
&\prod_{p=1}^n\ 
\biggl(\ 
\prod_{q=1}^{p-1}\,
\Rt_{1p}(u+c_q^{\,\prime}\com u+c_p^{\,\prime})\ot1
\biggr)
\ts F_p(u+c_p^{\,\prime})=
\\
&\prod_{p=1}^n\ 
\biggl(\ 
\prod_{q=1}^{p-1}\,
\Rt_{1p}(u+c_q\com u+c_p)\ot1
\biggr)
\ts F_p(u+c_p)
\cdot
X_{r,r+1}
\Tag{2.13}
\endalign
$$
in $\EndCN^{\ot n}\ot\Ug(u)$.
Combining the relations \(2.12)\ts,\ts\(2.13) we get
$$
F_{T^\prime}(u)=
\bigl(\ts X_{r,r+1}\ot1\ts\bigr)\,
\frac{F_T(u)}{1-d^{\ts2}}
\,\bigl(\ts X_{r,r+1}\ot1\ts\bigr)\,.
$$

We have the relation
$X^{2}=2d\ts X\hskip-1pt+\hskip-1pt1\hskip-1pt-d^{\ts2}$ in
the algebra $\EndCN^{\ot2}$.
The vectors $w_T$ and $w_{T^\prime}$ of the Young basis in the
$S_n$-module $W_\nu$ are orthogonal, so \(2.12) implies the equality
$Y_T\,X_{r,r+1}\,Y_T=0$. Therefore
$$
Y_T\,X_{r,r+1}^2\,Y_T=(1-d^{\ts2})\cdot Y_T\,.
$$
But the element $F_T(u)\in\EndCN^{\ot n}\ot\Ug(u)$ is divisible
by $Y_T\ot1$ on the right as well as on the left,
see the proof of Proposition 2.6\ts. Thus
$$
\gather
\,
(\tr^{\ts\ot n}\ot\id)\ts\bigl(F_{T^\prime}(u)\bigr)
=\ts
(\tr^{\ts\ot n}\ot\id)\,
\Bigl(\bigl(\ts X_{r,r+1}\ot1\ts\bigr)\,
\frac{F_T(u)}{1-d^{\ts2}}
\,\bigl(\ts X_{r,r+1}\ot1\ts\bigr)\Bigr)
\\
=\ts
(\tr^{\ts\ot n}\ot\id)\,
\Bigl(\bigl(\ts X_{r,r+1}^{\ts2}\ot1\ts\bigr)\,
\frac{F_T(u)}{1-d^{\ts2}}\,\ts\Bigr)
=\ts
(\tr^{\ts\ot n}\ot\id)\ts\bigl(F_T(u)\bigr)
\quad\square
\endgather
$$
\enddemo


\section{\bf\S\ts3. Leading terms of the element $Z_\nu(u)$}

Throughout this section $\nu=(\nu_1\ld\nu_N)$ will be any
partition of $n$ into not more than $N$ parts. However, we will always 
have $n=2m$. We fixed the basis $F_{-M,-M}\ld F_{-1,-1}\!$
in the Cartan subalgebra $\h\subset\g$.

Consider again the standard ascending filtration of the algebra $\Ug$
$$
\U_0(\g)\subset\U_1(\g)\subset\U_2(\g)\subset\ldots\subset\Ug\,.
$$
Here $\U_0(\g)=\CC\,$,$\U_1(\g)=\g$. By definition the subspace
$\U_n(\g)\subset\Ug$
consists of all the elements with degree not more than $n$. We will identify
the quotient space $\U_n(\g)/\U_{n-1}(\g)$ with the subspace in the
symmetric algebra $\Sg$ consisting of the homogeneous elements of degree $n$.
 
By \(2.6) we get $Z_\nu(u)\in\U_n(\g)\ot\CC(u)$.
The image of $Z_\nu(u)$ in
$$
\bigl(\ts\U_n(\g)/\U_{n-1}(\g)\bigr)\ot\CC(u)\,\subset\,\Sg\ot\CC(u)
$$
is a homogeneous polynomial in $F_{ij}\in\g$ of degree $n$
with the coefficients from $\CC(u)$. Due to Proposition 2.5 this
polynomial is invariant under the adjoint action of the group $G$
in $\Sg$. By the Chevalley theorem \text{[\ts D,\,Theorem 7.3.5\ts]}
this polynomial is uniquely determined by its image 
$$
f_\nu(x_1\ld x_M\ts|\ts u)\,\in\,\CC[\ts x_1\ld x_M]\ot\CC(u)
\Tag{3.1}
$$
with respect to the homomorphism $\eta:\Sg\to\CC[\ts x_1\ld x_M]$ defined by
the assignment $F_{ij}\mapsto0$ if $i\neq j$ or if $i=j=0$, and by
$$
F_{-M,-M}\mapsto x_1\ \ld\ F_{-1,-1}\mapsto x_M\,.
$$
Moreover, the image \(3.1) is a symmetric polynomial in
$x_1^{\ts2}\ld x_M^{\,2}$. Our present aim is
to determine the polynomial \(3.1) for any partition $\nu$ of $n=2m$.
In particular, we will describe the partitions $\nu$ where the polynomial
\(3.1) is not identically zero.

Denote by $\La_M$ the ring of symmetric polynomials in $x_1\ld x_M$
with complex coefficients.
For any partition $\rho=(\rho_1,\rho_2,\ts\ldots)$ put
$$
p_\rho(\ts x_1\ld x_M)\ts\,=\ts\,\prod_{k=1}^{\ts\ell(\rho)}\,\ts
\bigl(\ts x_1^{\ts\rho_k}+\ldots+x_M^{\ts\rho_k}\bigr)\,\in\,\La_M\,.
$$
As usual, here $\ell(\rho)$ is the number of non-zero parts in
the partition $\rho$.

We will use some
elementary facts from the representation theory of the 
symmetric group $S_{\ts2n}$. 
Consider the \text{hyperoctahedral} group $H_n$ as the subgroup in $S_{\ts2n}$
that centralizes the product of transpositions
$(1,n+1)\cdots(n,2n)\in S_{\ts2n}$.
Thus  $H_n=S_n\ltimes(\ZZ_2)^n$ where the subgroup
$S_n\subset S_{\ts2n}$ acts on $1\ld 2n$ by
simultaneous permutations of $1\ld n$ and $n+1\ld 2n$. Here
the subgroup $(\ZZ_2)^n\subset S_{\ts2n}$ is generated by 
the pairwise commuting transpositions
$(1,n+1)\ts\ld\ts (n,2n)$. Consider the one-dimensional representations
$\chi_+$ and $\chi_-$ of the group $H_n$ which are trivial on its subgroup
$S_n$ while
$
\chi_\pm:\,(s,n+s)\mapsto\pm1\ts
$
respectively.
Take the corresponding minimal idempotents in the group ring
$\CC\ts\!\cdot\!H_n$
$$
h_\pm\ts=\,\frac1{n\ts!\,2^n}\sum_{\si\in H_n}\chi_\pm(\si)\,\si\,.
$$
Note that the intersection of
$h_-(\CS_{\ts2n})h_+$ with $h_+(\CS_{\ts2n})h_-$
is zero.

\proclaim{\quad Proposition 3.1}
We can uniquely determine two linear maps
$$
\ch\hskip-1pt:\,h_-\ts(\CS_{\ts2n})\ts h_+\longrightarrow\La_M
\ \quad\text{and}\ \,\quad
\ch\hskip-1pt:\,h_+\ts(\CS_{\ts2n})\ts h_-\longrightarrow\La_M
$$
by setting
$$
\ch\ts(h_-\ts\si\ts h_+)=\ch\ts(h_+\ts\si\ts h_-)=
p_\rho(x_1^2\ld x_M^2)\cdot2^{\ts\ell(\rho)}
\hskip-40pt
$$
for any permutation $\si$ of $1\ld n$
with the cycle lengths $2\rho_1,2\rho_2\ts,\ts\ldots\,\,$.
\endproclaim

\demo{\quad Proof}
Any double coset of the subgroup $H_n$ in $S_{\ts2n}$ contains 
a permutation that acts on the numbers $n+1\ld 2n$ trivially.  
Moreover, all permutations $\si$
of $1\ld n$ with the same cycle lengths belong to the same double coset.
If any of these lengths is odd then
$h_-\ts\si\ts h_+=h_+\ts\si\ts h_-=0$.
Now for each partition $\rho=(\rho_1,\rho_2,\ldots)$ of $m$
choose a permutation $\si$ with
the cycle lengths $2\rho_1,2\rho_2\ts,\ts\ldots\,\,$.
All the corresponding elements $h_-\ts\si\ts h_+\in\CS_{\ts2n}$
are linearly independent. Therefore our definitions of 
two linear maps $\ch$ are self-consistent
\enddemos

We will call the two linear maps in Proposition 3.1 the
{\it characteristic maps,} see
\text{[\ts M\ts,\ts Section VII.2\ts].}
Now fix any standard tableau $T\in\Cal{T}_\nu$ and take the corresponding
minimal idempotent $\Phi_T\in\CS_n$.
Regard $\Phi_T$ as an element of the group ring
$\CS_{\ts2n}$ where the subgroup $S_n\subset S_{\ts2n}$
acts on the numbers $n+1\ld 2n$ trivially. Consider
the product
$$
\Psi_T(u)\,=\,
\prod_{p=1}^n\ \!
\Bigl(\ts1+\frac{(1,n+p)+\ldots+(p-1,n+p)}{2u+c_p}\ts\Bigr)
\cdot \Phi_T
\Tag{3.0}
$$
in $\CC(u)\!\ts\cdot\ts\!S_{\ts2n}$ where  
the (non-commuting) factors corresponding to the indices $p=1\ld n$
are as usual arranged from the left to the right.
Computation of the homogeneous
polynomial \(3.1) in $x_1\ld x_M$ hinges on the following observation.

\proclaim{\quad Proposition 3.2}
For $\g=\soN$ and $\g=\spN$
the polynomials \(3.1) coincide with the images in $\La_M\ot\CC(u)$
of $h_-\ts\Psi_T(u)\ts h_+$ and $h_+\ts\Psi_T(u)\ts h_-\!$
respectively under the characteristic maps.
\endproclaim

\demo{\quad Proof}
Take the permutational action of the group $S_{2n}$
in the space $(\CC^N)^{\ot\ts2n}$.  
Then the image in $\EndCN^{\ot\ts2n}(u)$ of $\Psi_T(u)$ is the product
$$
\prod_{p=1}^n\ \!
\Bigl(\ts1+\frac{P_{1,n+p}+\ldots+P_{p-1,n+p}}{2u+c_p}\ts\Bigr)
\cdot\bigl(\,Y_T\ot\id^{\ts\ot n}\ts\bigr)\,.
$$
Decompose this product with respect to the standard basis in
the space $\EndCN^{\ot\ts2n}\,$. We get the sum
$$
\sum_{i_1\ldots i_{2n}}\,\sum_{j_1\ldots j_{2n}}\ \ 
\psi^{\,i_1\ldots i_{2n}}_{j_1\ldots j_{2n}}(u)\cdot
E_{i_1j_1}\ot\cdots\ot E_{i_{2n}j_{2n}}
$$
where the coefficients 
are certain rational functions of $u$ valued in $\CC$.

We will put $\ep_i=\sgn i$ if $\g=\spN$ and set $\ep_i=1$ if $\g=\soN$.
Then $\ep_{ij}=\ep_i\ts\ep_j$ by definition. Denote
$$
I_{ij}(u)=\ep_i\cdot\bigl(\ts F_{j,-i}-(u+\eta)\,\de_{j,-i}\ts\bigr)
\in\Ug[\ts u\ts]\,,
\ \quad
J_{ij}=\ep_i\cdot\de_{i,-j}\,.
$$
By the definition of the element $F_T(u)\in\EndCN^{\ot n}\ot\Ug(u)$ we have 
$$
Z_\nu(u)=\!
\sum_{i_1\ldots i_{2n}}\ts\sum_{j_1\ldots j_{2n}}\  
\psi^{\,i_1\ldots i_{2n}}_{j_1\ldots j_{2n}}(u)
\cdot
I_{i_{n+1}i_1}(u)\ts J_{j_{n+1}j_1}
\cdots\,
I_{i_{2n}i_n}(u)\ts J_{j_{2n}j_n}
$$
where we employed the definition \(2.0) and the fact that
the elements $P,Q$ are obtained from each other
by applying to the second tensor factor of $\EndCN^{\ot2}$
the transposition with respect to $\langle\,\,,\,\rangle$.
This expression for $Z_\nu(u)$ shows that the polynomial
$f_\nu(x_1\ld x_M\ts|\ts u)$ equals the sum
$$
\sum_{i_1\ldots i_n}\,\sum_{j_1\ldots j_n}\ 
\psi^{\,i_1\ldots
i_n,-i_1\ldots-i_n}_{j_1\ldots j_n,-j_1\ldots-j_n}(u)
\cdot
\eta\ts(F_{i_1i_1})\,\ep_{i_1}\ep_{j_1}\cdots\,
\eta\ts(F_{i_ni_n})\,\ep_{i_n}\ep_{j_n}\,.
$$

The product
$
\eta\ts(F_{i_1i_1})\,\ep_{i_1}\ep_{j_1}
\cdots\,
\eta\ts(F_{i_ni_n})\,\ep_{i_n}\ep_{j_n}
$
is invariant under the permutations of 
the indices $i_1\ld i_n$ and of the indices $j_1\ld j_n$. 
For $\g=\soN$ it is also invariant under any substitution $j_p\mapsto-j_p$
with $p=1\ld n$
but changes the sign under the substitution $i_p\mapsto-i_p$.
Inversely, for $\g=\spN$ this product is
invariant under any substitution $i_p\mapsto-i_p$ but
changes the sign under the substitution $j_p\mapsto-j_p$.

Now take any
permutation of the first $n$ tensor factors in $(\CC^N)^{\ot\ts2n}$
with the cycle lengths $2\rho_1,2\rho_2\ts,\ts\ldots\,\,$
and decompose it in $\EndCN^{\ot\ts2n}$ as
$$
\sum_{i_1\ldots i_{2n}}\,\sum_{j_1\ldots j_{2n}}\  
\de^{\,i_1\ldots i_{2n}}_{\hskip.5pt j_1\ldots j_{2n}}\cdot
E_{i_1j_1}\ot\cdots\ot E_{i_{2n}j_{2n}}\!\!\!\!
$$ 
where each of the coefficients equals $0$ or $1$.
It remains to show that
$$
\sum_{i_1\ldots i_n}\,\sum_{j_1\ldots j_n}\ 
\de^{\,i_1\ldots i_n,-i_1\ldots-i_n}_{\hskip.5pt
j_1\ldots j_n,-j_1\ldots-j_n}
\cdot
\eta\ts(F_{i_1i_1})\,\ep_{i_1j_1}\cdots\,\eta\ts(F_{i_ni_n})\,\ep_{i_nj_n}
\Tag{3.2}
$$
then equals $p_\rho(x_1^2\ld x_M^2)\cdot2^{\ts\ell(\rho)}$
which is evident. Indeed, the latter expression and the sum \(3.2)
both are multiplicative with respect to the decomposition
of our permutation of the first $n$ tensor factors in $(\CC^N)^{\ot\ts2n}$
into the product of cycles. So we can assume that there is one single cycle
of length $n=2\ts m$. In this case
$$
\de^{\,i_1\ldots i_n,-i_1\ldots-i_n}_{\hskip.5pt
j_1\ldots j_n,-j_1\ldots-j_n}=
\cases
1&\text{\ if\ \ }i_1=j_1=\ldots=i_n=j_n\,,
\\
0&\text{\ otherwise}
\endcases
$$
and \(3.2) equals
$$
\sum_{i}\ \bigl(\,\eta\ts(F_{ii})\bigr)^n=
\,2\ts\bigl(\ts x_1^{\ts n}+\ldots+x_M^{\ts n}\bigr)\quad\square
\hskip-50pt
$$

Consider the two elements $h_-\ts\Psi_T(u)\ts h_+$
and $h_+\ts\Psi_T(u)\ts h_-$ of the ring
$\CC(u)\ts\!\cdot\ts\!S_{2n}$.
According to \text{Proposition 3.2\ts,} the first element
corresponds to the case $\g=\soN$ while the second corresponds to
$\g=\spN$.
We will evaluate the images of these two elements
under the corresponding characteristic maps by studying their actions 
in irreducible $S_{\ts2n}$-modules.

Let $\om$ be any partition of $2n$. The 
irreducible $S_{\ts2n}$-module $W_\om$ contains a non-zero vector $w_+$
such that $\si\cdot w_+=\chi_+(\si)\ts w_+$ for any $\si\in H_n$,
if and only if every row of the Young diagram of $\om$ has even length.
Then the vector $w_+\in W_\nu$ is unique up to a scalar multiplier,
and we will assume that $\langle\ts w_+\com w_+\rangle=1$.
The module $W_\om$ contains a non-zero vector $w_-$
with $\si\cdot w_-=\chi_-(\si)\ts w_-$ for any $\si\in H_n$,
if and only if every column of $\om$ has even length.
The vector $w_-\in W_\nu$ is then unique up to a scalar multiplier,
and we will assume that $\langle\ts w_-\com w_-\rangle=1$.
All these facts are well known, see for instance [\ts M, Section VII.2\ts].

Now suppose that
$\om=(\ts2\mu_1,2\mu_1,2\mu_2,2\mu_2,\ts\ldots)$ for a certain partition
$\mu=(\ts\mu_1,\mu_2,\ts\ldots)$ of $m$. We do not impose any restriction
on the number of parts in $\mu$ yet.
The partition $\om$ satisfies both conditions above, so
we have non-zero vectors $w_+,w_-\hskip-1pt\in W_\om\hskip-1pt$.
Let $b_1\ld b_m$ be the contents of the diagram $\mu$ ordered arbitrarily.

Let $\chi_\nu$ be the character of the irreducible $S_n$-module $W_\nu$,
take the element
$$
\Chi_\nu\,=\,\frac1{n!}
\sum_{\si\in S_n}\,
\chi_\nu(\si)\,\si
\,\in\CS_n\,.
$$
We will also regard $\Chi_\nu$ as an element of the group ring $\CS_{2n}$
by using the standard embedding $S_n\to S_{2n}$.
Where the double signs $\pm$ and $\mp$ appear in the next proposition, one
should simultaneously take only the upper signs or only the lower signs.
Recall that in this section $n=2m$.

\proclaim{\quad Proposition 3.3}
Action of the element
$h_\mp\ts\Psi_T(u)\ts h_\pm\in\CC(u)\cdot S_{2n}$
in the module $W_\om$ coincides with the action of
$h_\mp\Chi_\nu h_\pm\in\CS_{2n}$ times
$$
\frac
{(u+b_1)(u+b_1\pm1/2)\cdots(u+b_m)(u+b_m\pm1/2)}
{\hskip6pt(u+c_1/2)(u+c_2/2)\cdots(u+c_{n-1}/2)(u+c_n/2)}
\,\in\,\CC(u)\,.
\Tag{3.333}
$$
\endproclaim

\demo{\quad Proof}
For each $p=1\ld n$ consider the elements of the ring $\CS_{2n}$
$$
z_p^{\ts\prime}=\sum_{q=1}^{p-1}\ (n+q,n+p)
\quad\text{and}\quad
z_p^{\ts\prime\prime}\,=\,\sum_{q=1}^{p-1}\ (n+q,n+p)+(q,n+p)\,.
$$
The elements $z_1^{\ts\prime}\ld z_n^{\ts\prime}$
are the images in $\CS_{2n}$ of the Jucys-Murphy elements
$z_1\ld z_n\in\CS_n$
under the embedding
$$
S_n\to S_{\ts2n}:\ts(q,p)\mapsto(n+q,n+p)\,.
\Tag{3.22}
$$
In particular, the elements 
$z_1^{\ts\prime}\ld z_n^{\ts\prime}$ pairwise commute.
Note that the elements $z_1^{\ts\prime\prime}\ld z_n^{\ts\prime\prime}$
also pairwise commute.
The definition \(3.0) can be now rewritten as
$$
\Psi_T(u)\,=\,
\prod_{p=1}^n\ \!
\biggl(1+\frac{z_p^{\ts\prime\prime}-z_p^{\ts\prime}}{2u+c_p}\ts\biggr)
\cdot \Phi_T\,.
\Tag{3.3}
$$
But for any $p$ the element $z_p^{\ts\prime}$ commutes with
each of $z_{p+1}^{\ts\prime\prime}\ld z_n^{\ts\prime\prime}$.
On the other hand, due to Lemma 1.3 we have
the equalities in $\CS_{2n}$
$$
z_p^{\ts\prime}\ts\Phi_T\ts h_\pm=
\Phi_T\ts z_p\ts h_\pm=
c_p\cdot\Phi_T\ts h_\pm\,.
$$
Therefore \(3.3) implies the equality in  
the ring $\CC(u)\!\ts\cdot\ts\!S_{2n}$
$$
\Psi_T(u)\ts h_\pm=
\frac{(2u+z_1^{\ts\prime\prime})\cdots(2u+z_n^{\ts\prime\prime})}
{(2u+c_1)\cdots(2u+c_n)}
\,\ts\Phi_T\ts h_\pm\,.
\Tag{3.4}
$$

The standard chain of subgroups \(1.333) corresponds to the natural ordering
of the numbers $1\ld n$. Now consider the chain of subgroups
$$
S_1\subset S_2\subset\ldots\subset S_{2n-1}\subset S_{2n}
$$
corresponding the ordering
$n+1\ts,1\ts,n+2\ts,2\ts\ld 2n\ts,n$.
The elements $z_1^{\ts\prime\prime}\ld z_n^{\ts\prime\prime}\in\CS_{2n}$
are the Jucys-Murphy elements corresponding to the latter chain 
with the indices $1\ts,3\ts\ld 2n-1$.
Take the Young basis in $W_\nu$ corresponding to this chain of
subgroups in $S_{2n}$. The vectors $w_U$ of this basis are parametrized by
standard tableaux $U$ of shape $\om$ with the entries $1\ld 2n$.
But by [\ts BG,\,Theorem 3.4\ts] the vector $w_-\in W_\om$
is a linear combination of the vectors $w_U$ where $1\ts,3\ts\ld 2n-1$
occupy the first, third, $\ldots$ rows of the tableau $U$.
The collection of contents of the boxes in these rows is 
$2\ts b_1\ts,2\ts b_1+1\ts\ld 2\ts b_m\ts,2\ts b_m+1.$
By \text{Lemma 1.3} action of
$h_-\ts(2u+z_1^{\ts\prime\prime})\cdots(2u+z_n^{\ts\prime\prime})$ 
in $W_\om$ coincides with the action of
$$
(2u+2\ts b_1)(2u+2\ts b_1+1)\ldots(2u+2\ts b_m)(2u+2\ts b_m+1)\cdot h_-\,.
\hskip-10pt
$$
Similarly, the vector $w_+\in W_\om$
is a linear combination of the vectors $w_U$ where $1\ts,3\ts\ld 2n-1$
occupy the first, third, $\ldots$ columns of the tableau $U$.
The collection of contents of the boxes in these columns is 
$2\ts b_1\ts,2\ts b_1-1\ts\ld 2\ts b_m\ts,2\ts b_m-1.$
Again due to \text{Lemma 1.3} the action of
$h_+\ts(2u+z_1^{\ts\prime\prime})\cdots(2u+z_n^{\ts\prime\prime})$ 
in $W_\om$ coincides with the action of
$$
(2u+2\ts b_1)(2u+2\ts b_1-1)\ldots(2u+2\ts b_m)(2u+2\ts b_m-1)\cdot h_+\,.
\hskip-10pt
$$
Thus by \(3.4) the action of the element $h_\mp\ts\Psi_T(u)\ts h_\pm$
in the module $W_\nu$ coincides with the action of
$h_\mp\ts\Phi_T\ts h_\pm$ multiplied by the product \(3.333).

To complete the proof of Proposition 3.3 it remains to observe that
for any $\si\in S_n\subset S_{2n}$ we have
$h_\mp\ts\Phi_T\ts h_\pm=h_\mp\ts\si\ts\Phi_T\ts\si^{-1}h_\pm$. Therefore
$$
h_\mp\ts\Phi_T\ts h_\pm=
\frac1{n!}
\sum_{\si\in S_n}
h_\mp\ts\si\ts\Phi_T\ts\si^{-1}h_\pm=
h_\mp\Chi_\nu h_\pm
\quad\square
$$
\enddemo

Let us now formulate the main result of this section.
Consider again the ring $\La_N$ of symmetric polynomials
in the variables $y_1\ld y_N$. We assume that $x_1\ld x_M$
are independent of those $N$ variables.
Equip the vector space $\La_N$ 
with the standard inner product, so that
the Schur polynomials $s_\nu(\ts y_1\ld y_N)$ where $\nu$ runs
through the set of partitions
with not more than $N$ parts, constitute an orthonormal basis in $\La_N$.

Symmetric polynomial $s_\mu(\ts y_1^2\ld y_N^2)$ is the
{\it plethysm\/} of the Schur polynomial $s_\mu(\ts y_1\ld y_N)$ with the
power sum $y_1^2+\ldots+y_N^{\ts2}$. Expand
$$
s_\mu(\ts y_1^2\ld y_N^{\ts2})=\sum_\nu\ L_{\mu\nu}\,s_\nu(\ts y_1\ld y_N)
\Tag{3.5}
$$
in $\La_N$ with respect to the basis of Schur polynomials. The polynomials
$p_\rho(\ts y_1\ld y_N)$ form an othogonal basis in $\La_N$. If
$\rho_1+\rho_2+\ts\ldots=m$ and the number of permutations in $S_m$ with
the cycle lengths $\rho_1,\rho_2,\,\ldots$ is $m!\ts/\ts z_\rho$ then
the squared norm of $p_\rho(\ts y_1\ld y_N)$ is $z_\rho$. Further, then
$$
s_\mu(\ts y_1\ld y_N)\,=
\,\sum_\rho\ \chi_\mu^{\ts\rho}\,\ts
p_\rho(\ts y_1\ld y_N)/z_\rho
\Tag{3.51}
$$
where $\chi_\mu^{\ts\rho}$ denotes the value of the irreducible character
$\chi_\mu$
of $S_m$ on a permutation with the cycle lengths $\rho_1,\rho_2,\,\ldots$.
Therefore we have
$$
L_{\mu\nu}=\sum_\rho\ \chi_\mu^{\ts\rho}\,\chi_\nu^{\ts2\rho}/\ts z_\rho\,.
\Tag{3.6}
$$
As usual, we denote $2\rho=(\ts2\rho_1,2\rho_2\ts,\ts\ldots)$.
Note that $z_{\ts2\rho}=2^{\ts\ell(\rho)}z_\rho$ then.

Combinatorial description of the
coefficients $L_{\mu\nu}$ in the expansion \(3.5)
has been provided in [\ts L,\,Section 5\ts].
Another description
of these coefficients is given by [CL,\,Theorem 5.3\ts].
In particular, if $L_{\mu\nu}\neq0$ then the Young diagram of $\nu$
can be split into horizontal and vertical blocks of two boxes each.
These blocks are called {\it dominoes\/}, see [BG].
\Par
We put
$\eta=\frac12$ if $\g=\soN$ and put $\eta=-\frac12$ if $\g=\spN$.
Recall that $c_1\ld c_n$ are the contents of a standard tableau $T$ of shape
$\nu$. The contents $b_1\ld b_m$ of the boxes of $\mu$
have been ordered arbitrarily.

\proclaim{\quad Theorem 3.4}
The polynomial $f_\nu(x_1\ld x_M\ts|\ts u)$ equals the sum over
all partitions $\mu$ of $m$ into not more than $M$ parts, of the products
$$
\frac
{(u+b_1)(u+b_1+\eta)\cdots(u+b_m)(u+b_m+\eta)\hskip6pt}
{(u+c_1/2)(u+c_2/2)\cdots(u+c_{n-1}/2)(u+c_n/2)}
\,
L_{\mu\nu}
\,
s_\mu(\ts x_1^2\ld x_M^{\ts2})
\,.
$$
\endproclaim

\demo{\quad Proof}
In this proof the upper signs in $\pm$ and $\mp$ correspond to
the case $\g=\soN$ while the lower signs correspond to  $\g=\spN$.
Initially let $\mu$ run through the set of all partitions of $m$,
without any restriction on the number of parts. The elements
$$
\Gm\,=\,\frac{\dim W_\om}{(2n)!}
\sum_{\tau\in S_{2n}}\,
\langle\,w_\mp\com\tau\!\cdot\!\hskip.5pt w_\pm\,\rangle_\om
\,h_\mp\ts\tau\ts h_\pm
\,\in\,\CS_{2n}
\Tag{3.61}
$$
form a basis in the vector space $h_\mp\ts(\CS_{2n})\ts h_\pm$. 
Let us expand
$$
h_\mp\ts\Psi_T(u)\ts h_\pm=\sum_\mu\ 
f_{\mu\nu}(u)\,\Gm
\nopagebreak
\Tag{3.62}
$$
with respect to this basis and compute the coefficients
$f_{\mu\nu}(u)\in\CC(u)$.
The element $h_\mp\ts\tau\ts h_\pm$ acts in the $S_{2n}$-module $W_\om$
as the linear operator
$\langle\,\tau\!\cdot\!\hskip.5pt w_\pm\com w_\mp\ts\rangle_\om\ts E$
where
$E:\ts w\mapsto\langle\,w\com w_\pm\,\rangle_\om\,w_\mp$
for any vector $w\in W_\om$.

The element $\Gm$ acts as the operator $E$ in the module $W_\om$ and
vanishes in any other irreducible $S_{2n}$-module. 
Denote by $d_{\mu\nu}(u)$ the rational function \(3.333).
By Proposition 3.3 and by the definition of $\Chi_\nu$ 
$$
f_{\mu\nu}(u)
\,=\,\frac{d_{\mu\nu}(u)}{n!}
\sum_{\si\in S_n}\,
\chi_\nu(\si)\,
\langle\,\si\!\cdot\!\hskip.5pt w_\pm\com w_\mp\ts\rangle_\om\,.
\Tag{3.7}
$$
Here the factor
$\langle\,\si\!\cdot\!\hskip.5pt w_\pm\com w_\mp\ts\rangle_\om$
may be non-zero only if the permutation $\si$ has the cycle lengths
$2\rho_1,2\rho_2,\,\ldots$ for some partition $\rho$ of $m$. Then 
$$
\langle\,\si\!\cdot\!\hskip.5pt w_\pm\com w_\mp\ts\rangle_\om=
I_\mu\cdot2^{\ts\ell(\rho)}\ts\chi_\mu^{\ts\rho}
\Tag{3.71}
$$
where $I_\mu$ depends only on the choice of the vectors $w_+\com w_-\in W_\om$
and
$$
|\ts I_\mu|^2=\frac{(2n)!}{\dim W_\om\cdot(\ts2^n\ts n!\ts)^2}\,.
\Tag{3.72}
$$
This result was independently obtained by Ivanov [\ts I,\,Theorem 3.9\ts]
and Rains [\ts R,\,Corollary 7.6\ts].
Using \(3.6) and \(3.7) along with this result,
$$
f_{\mu\nu}(u)\,=\,
d_{\mu\nu}(u)\ts
I_\mu\cdot
\sum_\rho\ \chi_\mu^{\ts\rho}\,\chi_\nu^{\ts2\rho}/\ts z_\rho
\,=\,
d_{\mu\nu}(u)\ts
I_\mu\ts L_{\mu\nu}\,.
\Tag{3.8}
$$

To complete the proof of Theorem 3.4
it now remains to apply the characteristic
map to each side of the equality \(3.62). By Proposition 3.2
on the left-hand side we get the polynomial
$f_\nu(x_1\ld x_M\ts|\ts u)$. 
There are exactly
$(2^n\ts n!)^2/\ts(4^{\ts\ell(\rho)}z_\rho)$ elements
in the double coset of the subgroup $H_n$ in $S_{2n}$
containing the permutation of $1\ld n$ with the cycle lenghts
$2\rho_1,2\rho_2,\,\ldots\,$.
By the definition \(3.61) and again by \text{\(3.71),\(3.72)}
$$
\ch\ts(\ts\Gm)
\,=\,
I_\mu^{\ts-1}\cdot
\sum_\rho\ 
\chi_\mu^{\ts\rho}\ts\,p_\rho(x_1^2\ld x_M^2)/z_\rho
\,=\,
I_\mu^{\ts-1}\cdot
s_\mu(\ts x_1^2\ld x_M^{\ts2})\,,
\nopagebreak
$$
here we have also used Proposition 3.1 and the classical expansion~\(3.51).
Thus the expression \(3.8) for the coefficient in \(3.62) shows that 
$$
\align
f_\nu(x_1\ld x_M\ts|\ts u)
\,=\,
&\sum_\mu\ d_{\mu\nu}(u)\,I_\mu\,L_{\mu\nu}\cdot\ch\ts(\ts\Gm)
\,=\,
\\
\sum_\mu\ d_{\mu\nu}(u)\,&L_{\mu\nu}\,s_\mu(\ts x_1^2\ld x_M^{\ts2})\,.
\endalign
$$
The latter sum can be restricted to the partitions $\mu$ with not more than
$M$ parts, since for the other partitions we have 
$s_\mu(\ts x_1^2\ld x_M^{\ts2})=0$
\enddemos

\proclaim{\quad Corollary 3.5}
If the polynomial \(3.1) corresponding to $\nu$
is not identically zero, then the
Young diagram of $\nu$ splits into dominoes.
\endproclaim

One can reformulate Theorem 3.4 as follows,
cf.\ [OO2\ts,\ts Theorem 1.2]\ts. Let us denote by
$b_\mu(u)$ and $c_\nu(u)$
the numerator and the denominator of the fraction in \(3.333).
The upper signs in $b_\mu(u)$ correspond to 
$\g=\soN$ while the lower signs correspond to $\g=\spN$.

\proclaim{\quad Corollary 3.6}
For any fixed positive integers $m$ and $N$ we have
$$
\gather
\sum_\nu\,\,c_\nu(u)\,f_\nu(x_1\ld x_M)\,s_\nu(y_1\ld y_N)=
\\
\sum_\mu\,\,b_\mu(u)\,
s_\mu(x_1^{\ts2}\ld x_M^{\ts2})\,s_\mu(y_1^{\ts2}\ld y_N^{\ts2})\,.
\Tag{3.10}
\endgather
$$
where $\nu$ and $\mu$ range respectively over all partitions 
of $n=2m$ with at most $N$ parts and all partitions
of $m$ with at most $M=[N/2\ts]$ parts.
\endproclaim

\demo{\quad Proof}
By Theorem 3.4 for any partition $\nu$ of $n=2m$ into not more than $N$
parts the product $c_\nu(u)\,f_\nu(x_1\ld x_M)\,s_\nu(y_1\ld y_N)$ equals
$$
\sum_\mu\,\,b_\mu(u)\,L_{\mu\nu}\,
s_\mu(x_1^{\ts2}\ld x_M^{\ts2})\,s_\nu(y_1\ld y_N)\,.\
$$
Taking here the sum over $\nu$ we obtain \(3.10) by the definition \(3.5)
\enddemos

We will complete this article with the following two examples.
Firstly, let us put
$\mu=(m\ts,0\ts\ld0)$ and $\nu=(2m\ts,0\ts\ld0)$.
Then the element $B_\mu\in\Zg$ described in Section 2, coincides with
the value of $Z_\nu(u)$ at $u=-\ts m-\frac12$ for $\g=\soN$ and with
the value of $(u+m-\frac12)/(u-\frac12)\ts\cdot\ts Z_\nu(u)$
at the point $u=-\ts m\!+\!\frac12$ for $\g=\spN$;
see [\ts MN,Theorem 3.3\ts]\ts.

Secondly, put
$\mu=(1\ts\ld1\ts,0\ts\ld0)$ and $\nu=(1\ts\ld1\ts,0\ts\ld0)$
where the part $1$ appears $m$ and $2m$ times respectively. Then
the element $(-1)^m B_\mu\in\Zg$ coincides with
the value of $Z_\nu(u)$ at $u=m+\frac12$ for $\g=\spN$ and with
the value of $(u-m+\frac12)/(u+\frac12)\cdot Z_\nu(u)$
at the point $u=m-\frac12$ for $\g=\soN$;
see [\ts MN,Theorem 3.4\ts]\ts.

It would be interesting to establish a connection between our functions
$Z_\nu(u)$ and the elements $B_\mu\in\Zg$ with the general partitions $\mu$.


\vbox{
\bigskip\centerline{\bf Acknowledgements}
\section{\,}\kern-20pt

I am grateful to A.\,Lascoux, B.\,Leclerc and J.\ts-Y.\,Thibon
for valuable remarks. I am especially indebted to G.\,Olshanski.
Discussions with him of the results [OO2] have inspired
the present work. Financial support from the EPSRC and from the EC
under the grant FMRX-CT97-0100 is gratefully acknowledged.
}


\bigskip
\centerline{\bf References}
\section{\,}\kern-20pt

\itemitem{[BG]}
{N.\,Bergeron\,and\,A.\,Garsia},
{\!Zonal\,polynomials and\,domino\,tableaux},
{Discrete Math.},
{\bf 99}
(1992),
3--15.

\itemitem{[C]}
{A.Capelli},
{Sur les op\'erations dans la th\'eorie des formes alg\'ebriques},
{Math.\,Ann.},
{\bf 37}
(1890),
1--37.

\itemitem{[C1]}
{I.\,Cherednik},
{Factorized particles on the half-line and root systems},
{Theor.\,Math.\,Phys.},
{\bf 61}
(1984),
977--983.

\itemitem{[C2]}
{I.\,Cherednik},
{On special bases of irreducible finite-dimensional representations
of the degenerate affine Hecke algebra},
{Funct.\,Analysis Appl.},
{\bf 20}
(1986),
87--89.

\itemitem{[CL]}
{C.\,Carr\'e and B.\,Leclerc},
{Splitting the square of a Schur function into
its symmetric and antisymmetric parts\/},
{J. Algebraic Combin.},
{\bf 4}
(1995),
201--231.

\itemitem{[D]}
{J.\,Dixmier},
{\!\lq\lq Alg\`ebres enveloppantes\rq\rq\!},
{Gauthier-Villars,\,Paris},\,
1974.

\itemitem{[I]}
{V.\,Ivanov},
{Bispherical functions on the symmetric group,
associated to the hyperoctahedral subgroup},
to appear in J.\,Math.\,Sci.

\itemitem{[J]}
{A.\,Jucys},
{Symmetric polynomials and the centre of the symmetric group ring\/},
{Rep.\,Math.\,Phys.},
{\bf 5}
(1974),
107--112.

\itemitem{[JKMO]}
{M.\,Jimbo, A.\,Kuniba,\,T.\,Miwa and M.\,Okado},
{The $\!A_n^{(1)}\!\!$ face \text{models}},
{Comm.\,Math.\,Phys.},
{\bf 119}
(1988),
543--565.

\itemitem{[L]}
{D.\,Littlewood},
{Modular representations of symmetric groups\/},
{Proc. Royal Soc.},
{\bf A\,209}
(1951),
333--353.

\itemitem{[M]}
{I.\,Macdonald},
{\!\lq\lq Symmetric Functions and \hskip-.3pt Hall Polynomials\rq\rq\!},
\!Claren\-don Press, Oxford, 1995.

\itemitem{[M1]}  
{A.\,Molev},
{Sklyanin determinant, Laplace operators, and characteristic identities
for classical Lie algebras},
{J.\,Math.\,Phys.},
{\bf 36}
(1995),
923--943.

\itemitem{[M2]}
{A.\,Molev},
{Factorial supersymmetric Schur functions and super Ca\-pel\-li identities},
{\lq\lq\ts Kirillov's Seminar on Representation Theory\ts\rq\rq},
edited by G.\,Olshanski,
Amer.\,Math.\,Soc.\,Translations,
{\bf 181}
(1998),
109--137.

\itemitem{[MN]}
{A.\,Molev and M.\,Nazarov},
{Capelli identities for classical Lie algebras},
to appear in Math.\,Ann.

\itemitem{[MNO]}
{A.\,Molev, M.\,Nazarov and G.\,Olshanski},
{Yangians and classical Lie algebras},
Russian Math.\,Surveys,
{\bf 51}
(1996),
205--282.

\itemitem{[N1]}
{M.\,Nazarov},
{Quantum Berezinian and the classical Capelli identity},
{Lett.\,Math.\,Phys.},
{\bf 21}
(1991),
123--131.

\itemitem{[N2]}
{M.\,Nazarov},
{Yangians and Capelli identities},
{\lq\lq\ts Kirillov's Seminar on Representation Theory\ts\rq\rq},
edited by G.\,Olshanski,
Amer.\,Math.\,Soc. Translations,
{\bf 181}
(1998),
139--163.

\itemitem{[O]}
{A.\,Okounkov},
{Quantum immanants and higher Capelli identities},
{Transformation Groups},
{\bf 1}
(1996),
99--126.

\itemitem{[O1]}
{G.\,Olshanski},
{Twisted Yangians and infinite-dimensional classical Lie algebras},
\lq\lq\ts Quantum Groups\ts\rq\rq\,,
edited by P.\,Kulish,
Lecture Notes in Math.,
{\bf 1510}
(1992),
103--120.

\itemitem{[O2]}
{G.\,Olshanski},
{Generalized symmetrization in enveloping algebras},
{Transformation Groups},
{\bf 2}
(1997),
197--213.

\itemitem{[OO1]}
A.\,Okounkov and G.\,Olshanski,
{Shifted Schur functions},
St.\,Peters\-burg Math.\,J.,
{\bf 9}
(1998),
239--300.

\itemitem{[OO2]}
A.\,Okounkov and G.\,Olshanski,
{Shifted Schur functions II.
The binomial formula for characters of classical groups and its applications},
{\lq\lq\ts Kirillov's Seminar on Representation Theory\ts\rq\rq},
edited by G.\,Olshanski,
Amer.\,Math.\,Soc.\,Translations,
{\bf 181}
(1998),
245--271.

\itemitem{[R]}
E.\,Rains,
{\lq\lq\ts Attack of the Zonal Polynomials\ts\rq\rq},
Harvard University preprint, 1995.

\itemitem{[S]}
{S.\,Sahi},
{The spectrum of certain invariant differential operators associated
to a Hermitian symmetric space,}
{\lq\lq\ts Lie Theory and Geometry\ts\rq\rq},
edited by J.-L.\,Brylinski, R.\,Brylinski, V.\,Guillemin and V.\,Kac,
Progress in Math.,
{\bf 123}
(1994),
569--576.

\itemitem{[S1]}
E.\,Sklyanin,  
{Boundary conditions for integrable equations,}
{Funct. Analysis Appl.},
{\bf 21}  
(1987),
164--166.
  
\itemitem{[S2]}
E.\,Sklyanin,
{Boundary conditions for integrable quantum systems,}
J.\,Phys.,
{\bf A\,21}
(1988),
2375--2389.  

\itemitem{[W]}
{H.\,Weyl},
\text{\!\lq\lq
Classical Groups, their Invariants and Representations\rq\rq},
Princeton University Press, Princeton, 1946.

\itemitem{[Y1]}
{A.\,Young,}
{On quantitative substitutional analysis I and II\ts},
{Proc. London Math.\,Soc.},
{\bf 33}
(1901),
97--146
and
{\bf 34}
(1902),
361--397.

\itemitem{[Y2]}
{A.\,Young,}
{On quantitative substitutional analysis VI},
{Proc.\,London Math.\,Soc.},
{\bf 34}
(1932),
196--230.


\vbox{
\bigskip
\line{\it Department of Mathematics\hfill}
\line{\it University of York\hfill}
\line{\it York YO1 5DD, England\hfill}
\vskip-20pt
}


\bye